\documentclass[11pt,reqno]{article}

\usepackage{amssymb,mathtools,amsfonts,amsmath}
\usepackage{amsthm}
\usepackage{graphicx}
\usepackage{xcolor}
\usepackage{microtype}
\usepackage[colorlinks, linkcolor=black, anchorcolor=black,citecolor=black]{hyperref}

\usepackage[a4paper, margin=25mm]{geometry}

\allowdisplaybreaks
\numberwithin{equation}{section}

\theoremstyle{plain}
\newtheorem{theorem}{Theorem}[section]
\newtheorem{proposition}{Proposition}[section]
\newtheorem{lemma}{Lemma}[section]
\newtheorem{corollary}{Corollary}[section]
\theoremstyle{definition} 
\newtheorem{definition}{Definition}[section]
\newtheorem{assumption}{Assumption}[section]
\theoremstyle{remark}
\newtheorem{remark}{Remark}[section]
 
\newcommand{\R}{\mathbb{R}}

\newcommand{\csh}{c_{\rm sh}}
\newcommand{\cstar}{c_*}
\newcommand{\av}[1]{\left\langle #1\right\rangle}

\begin{document}
\title{ Asymptotic spreading  of  Fisher-KPP equations in shifting environments with time-periodic limits}
\author{Zhucheng Jin and Zhenkun Wang\thanks{ Corresponding author. \\ Email: jinzc@ustc.edu.cn (Z. Jin); zhenkunwang@ustc.edu.cn (Z. Wang)}\\
\small  School of Mathematical Sciences, University of Science and Technology of China,  
 Hefei 230026, P.R. China}
\date{}
\maketitle

\begin{abstract}
 In this paper, we study the rightward spreading  speed of solutions  for Fisher–KPP equations in shifting environments. 
 The linear growth rate has time periodic
 limits at the two ends of the moving frame. 
 Each limit has a positive average and may be  negative during   a period. 
 When the growth rate lies between the two limits at each
 time, we show that the spreading speed depends only on their averages and the asymptotic shift speed. 
 For intermediate shift speeds, the formulas give nonlocal pulling if the right
 end has the larger average and speed locking if the left end has the larger average. The interior
 growth profile may be nonmonotone in space and aperiodic in time. 
 We prove the formulas by constructing upper and
 lower solutions for a stepwise profile and then comparing with shifted stepwise problems.
\end{abstract}



\section{Introduction}
Seasonal variations in population growth and the geographic shift of habitat boundaries often occur together during species invasion.  
When the spatial distribution of resources or climatic conditions changes in a landscape,   a species may experience alternating  periods of growth and decline. Reaction-diffusion models in shifting environments have been widely studied to understand how habitat motion affects population persistence and invasion, since  \cite{BerestyckiDiekmannNagelkerkeZegeling2009, PotapovLewis2004}. When the environment also fluctuates periodically in time, 
a natural question is how to   determine the invasion speed.   
Because growth  ahead of and behind the moving interface  follow different seasonal patterns, we ask when the  temporal  averages of the limiting growth rates suffice to  give the correct invasion speed and whether growth within the transition zone can change the speed even under  the  fixed  limiting growth rates and the shift speed.

We study the Cauchy problem
\begin{equation}\label{Pb}
\begin{cases}
 u_t=u_{xx}+u f(t,x-A(t),u),&t>0, x\in\R, \\
 u(0,x)=u_0(x),&x\in\R,
\end{cases}
\end{equation}
where $u(t,x)$  is the density of  population, $f$  is  the  per capita 
growth rate and    $u_0$ is nonnegative, nontrivial   and compactly supported.
The function $A(t)$ describes the moving  of the environment.   The reaction term is of the Fisher-KPP type. 
  Set
$r(t,z)=f(t,z,0)$.  Assume that $r$ has two continuous
$T$-periodic spatial limits $r_-(t)$ and $r_+(t)$, uniformly in time. 
Their  averages   are  required to be positive, namely, 
$$
 R_\pm: =\frac1T\int_0^T r_\pm(s)\,ds>0.
$$
The instantaneous growth rate can  be negative during  the period.

\vspace{1ex}
For a classical Fisher-KPP equation in a homogeneous environment
$$u_t=  u_{xx} +ru(1-u),    $$
the asymptotic spreading speed for the Cauchy problem  is proved to be $2\sqrt{r}$, see \cite{Fisher1937,KPP1937,AronsonWeinberger1978}.  For time periodic equations, the spreading obtained by replacing $r$ with its temporal average. General propagation theories for time periodic and  general time heterogeneity of  KPP equations were developed in \cite{liang2006spreading,LiangZhao2007,NadinRossi2012,shen2010variational, WLiZ2026}.

Introducing a moving spatial transition adds a further velocity scale.    Such models were  first considered in \cite{BerestyckiDiekmannNagelkerkeZegeling2009, PotapovLewis2004} to describe the climate changing.  
 Beyond the biological motivation, the shifting environment appears in a triangular system that is partially decoupled, see \cite{FangLouWu2016,HolzerScheel2014}.
Forced waves and their stability in autonomous shifting environments have also been extensively studied in \cite{BerestyckiFang2018, NiuFang2025,  WangZhao2019, Zhao2026}.
When the interface of a shifting habitat  separates a favorable region from an unfavorable one, or when the favorable habitat forms a moving patch of bounded width, the  sharp criteria for 
species persistence, critical patch sizes, and asymptotic spreading speeds  were  established in \cite{LiBewickShangFagan2014, HuLi2015}.   Propagation phenomena of reaction-diffusion systems  in  shifting environments were   studied in \cite{DongShangLi2021,GuoHamel2023,Jin2026, WPOu2022, YiZhao2020}.

 For time-periodic shifting environments connecting favorable and unfavorable limiting habitats, Fang, Peng and Zhao \cite{FangPengZhao2021}  established the existence, uniqueness, and attractivity of time-periodic forced waves together with spreading properties.  Some dynamical methodology were developed to investigate the propagation phenomena for monotone and non-monotone evolution systems without spatial translation invariance, see \cite{YiZhao2020,YiZhao2025, HouWangZhao2025}.  Recently, the spreading properties and forced traveling waves in general time heterogeneous shifting environments was further studied by \cite{ZhangZhao2025}.

For the shifting  environments connecting two favorable ends,  
 the spreading speed may be not locally determined under certain regime,   see \cite{HuShangLi2020}. 
   The  precise nonlocal pulled speed  were given by \cite{LamYu2022} through developing a Hamilton-Jacobi approach.  An alternative methodology,   given in Chapter~6 of \cite{LamLou2022}, is  based on gluing elementary upper and lower solutions.  The \emph{nonlocal pulling} phenomena means that the spreading speed can  be enhanced by the presence of a good habitat at far distance.   This phenomenon also observed  in competition systems  \cite{GLam2019} and  a model for describing  spread of a trait in a population, introduced in \cite{V2014}.   A probabilistic method was used in \cite{BovH2023} to analyze the nonlocal pulling phenomena  in the genetic model in \cite{V2014}. 
     Recently, the influence of the moving interior profile  was further studied by \cite{GilettiGirardinMatano2024} and \cite{LamNadinYu2025}. 
   The  nonmonotone piecewise constant profile was considered by  \cite{GilettiGirardinMatano2024} through constructing  suitable super and subsolutions.  The   continuous case was investigated in \cite{LamNadinYu2025} by using a Hamilton–Jacobi approach with a   dynamic junction condition.  
  For more references, we refer the reader to surveys \cite{WL2022,LamWuYu2026}.

   \vspace{1ex}

In this work, we  focus on   that the linear growth  rate lies between its two spatial limits:
$$
 r_-(t)\leq r(t,z)\leq r_+(t)  \qquad \forall t\in \R
 \qquad\text{or}\qquad
 r_+(t)\leq r(t,z)\leq r_-(t)   \qquad \forall t\in \R. 
$$
   The spatial profile of $f$ may be nonmonotone. For each given  $z\in \R$ and $u>0$,  the growth rate function   $ f(t, z, u)$ can vary nonperiodically in time.    We assume only that $A$ is continuous and  $A(t)/t\to\csh>0$.    The function $A(t)-\csh t$ may be unbounded. 
 The precise assumptions  of $f$  are given in  Section~\ref{Sec-results}.
 
We establish explicit rightward spreading speed formulas when both limiting habitats are favorable on average and their growth rates vary periodically in time. The spreading speed  depends only on $R_-$, $R_+$ and $\csh$.   If $ r_-(t)\leq r(t,z)\leq r_+(t)$, the spreading speed is  nonincreasing as $\csh$ increases. Moreover, in the middle region,  the  nonlocally pulling phenomena occurs  with  $2\sqrt{R_-}<\cstar<\csh$.
Under the  decreasing ends  order case, 
the spreading   speed equals  $\csh$
when $\csh \in (2\sqrt{R_+}, 2\sqrt{R_-})$. So  in the middle region,  the
asymptotic invasion speed follows the shifting speed of the environment. This is the  speed
locking phenomenon.

These formulas identify a class of seasonal shifting environments
in which the two average of  growth rates and the asymptotic shift speed
are sufficient to determine spreading. 
Within this class, temporal amplitudes, the interior profile, and nonlinear
reaction  may change without influencing the spreading speed, provided
the hypotheses and  parameters $R_{\pm}$ and $\csh$ are preserved.

The spreading speed formulas is established through explicit
upper and lower solutions.   We investigate stepwise reactions first and then extend the  construction to growth rates  lying between their two spatial limits.   A change of time variable also gives an exact speed formula for spatially homogeneous
time-periodic diffusion. For shifts without a limiting velocity,
we derive the bounds of  spreading speed  from the lower and upper limits of
$A(t)/t$.

The remainder of the paper is organized as follows. Section~\ref{Sec-results} states the assumptions and main results. Section~\ref{Sec-step} analyzes the auxiliary stepwise problem, and Section~\ref{Sec-general} completes the proof of the main theorems. Section~\ref{Sec-consequences} extends the main results to time-periodic diffusion and  shifts for which $A(t)/t$ need not converge.

\section{Assumptions and main results}\label{Sec-results}
For a continuous $T$-periodic function $b$, we set 
$$
 \av b: =\frac1T\int_0^T b(s)\,ds.
$$
\begin{assumption}\label{ASS-f}
Let $f:\R^2\times[0,\infty)\to\R$ be continuous. Set $r(t,z)= f(t,z,0)$.  Assume that  there are constants
$M,B,L,K,\sigma_0>0$ such that, for all $t,z\in\R$ and $u,v\in[0,M]$,
\begin{align}\label{eq-dissp}
 |f(t,z,u)|&\leq B,
 &|f(t,z,u)-f(t,z,v)|\leq L|u-v|,\nonumber\\
 f(t,z,M)&\leq-\sigma_0,     \quad    & r(t,z)-Ku\leq f(t,z,u)\leq r(t,z) ,    
\end{align}
\end{assumption}

 Here $f(t,z,u)$ is not assumed to be monotone in $u$ or periodic in $t$.  Although the lower bound in the KPP condition follows from the Lipschitz condition with $K=L$, we keep $K$  as the constant used in the subsolution estimates. 
 
\begin{assumption}\label{ASS-limits}
There exist continuous $T$-periodic functions $r_-,r_+$ such that
\begin{equation}\label{eq-limitsf}
 \lim_{z\to\pm\infty}\sup_{t\in\R}|r(t,z)-r_\pm(t)|=0.
\end{equation}
Assume $R_\pm:=\av{r_\pm}>0. $
One of the following alternatives holds for every $(t,z)\in\R^2$:
\begin{align*}
 (H_\uparrow)&\qquad r_-(t)\leq r(t,z)\leq r_+(t),
 &0<R_-<R_+,\\
 (H_\downarrow)&\qquad r_+(t)\leq r(t,z)\leq r_-(t),
 &0<R_+<R_-.
\end{align*}
\end{assumption}

\begin{remark}
 Condition \eqref{eq-limitsf} relates to  the spatial limits of $f(t,z,u)$ only at $u=0$, without requiring spatial limits or a pointwise ordering for $u>0$. We can observe that $r(t,z)$ and $r_{\pm}(t)$ are bounded.  The function $r_-(t)$ and $r_+(t)$ may change sign. 
\end{remark}

\begin{assumption}\label{ASS-A}
The function $A\in C([0,\infty))$ satisfies
\begin{equation}\label{eq-A}
 \lim_{t\to\infty}\frac{A(t)}t=\csh>0.
\end{equation}
\end{assumption}

We next define the rightward spreading speed.
\begin{definition}\label{DEF-speed}
A number $c_*>0$ is the rightward spreading speed of \eqref{Pb} if
\begin{align*}
 &\lim_{t\to\infty}\sup_{x\geq ct}u(t,x)=0
       &&\forall  c>c_*,\\
 &\liminf_{t\to\infty}\inf_{0\leq x\leq ct}u(t,x)>0
       && \forall  0<c<c_* .
\end{align*}
Define   the lower and upper spreading speeds by
\begin{align}
 \underline c(u) 
 &:=\sup\left\{c>0:
   \liminf_{t\to\infty}\inf_{0\leq x\leq ct}u(t,x)>0\right\},
   \label{eq-lower-speed}\\
 \overline c (u)
 &:=\inf\left\{c>0:
   \lim_{t\to\infty}\sup_{x\geq ct}u(t,x)=0\right\}.
   \label{eq-upper-speed}
\end{align}
\end{definition}



We now state the main results.
\begin{theorem}\label{Thm-inc}
Let Assumptions \ref{ASS-f}, \ref{ASS-limits} and \ref{ASS-A}, with alternative
  $(H_\uparrow)$,  be satisfied.  Assume
$u_0\in C (\R,[0,M])$, $u_0\not\equiv0$ and compactly supported.
Then the solution of \eqref{Pb} has the rightward spreading speed
\begin{equation*}
 \cstar=C_{\uparrow}(\csh;R_-,R_+):=
 \begin{cases}
  2\sqrt{R_+},
    &0<\csh\leq2\sqrt{R_+},\\[2mm]
  \displaystyle \lambda_*+\frac{R_-}{\lambda_*},
    &2\sqrt{R_+}<\csh<
       2\bigl(\sqrt{R_+-R_-}+\sqrt{R_-}\bigr),\\[2mm]
  2\sqrt{R_-},
    &\csh\geq
       2\bigl(\sqrt{R_+-R_-}+\sqrt{R_-}\bigr),
 \end{cases}
\end{equation*}
where
$$
 \lambda_*=\frac{\csh}{2}-\sqrt{R_+-R_-}.
$$
\end{theorem}

\begin{remark}
If $2\sqrt{R_+}<\csh<
2\bigl(\sqrt{R_+-R_-}+\sqrt{R_-}\bigr)$, then
$$2\sqrt{R_-}<\cstar<\csh.  $$
In this interval, $\cstar$ depends on both $R_-$ and
$R_+$ and  is  larger than  $2\sqrt{R_-}$. This corresponds to the nonlocal  pulling phenomena, see also \cite{HolzerScheel2014,GilettiGirardinMatano2024,LamYu2022,LamNadinYu2025}.
\end{remark}

\begin{theorem}\label{Thm-dec}
Let Assumptions \ref{ASS-f}, \ref{ASS-limits} and \ref{ASS-A}, with
alternative $(H_\downarrow)$, be satisfied.  Assume
$u_0\in C (\R,[0,M])$, $u_0\not\equiv0$ and compactly supported.
Then the rightward spreading speed of \eqref{Pb} is
$$
 \cstar=C_{\downarrow}(\csh;R_-,R_+):=
 \begin{cases}
  2\sqrt{R_+},&0<\csh\leq2\sqrt{R_+},\\
  \csh,&2\sqrt{R_+}<\csh<2\sqrt{R_-},\\
  2\sqrt{R_-},&\csh\geq2\sqrt{R_-}.
 \end{cases}
$$
\end{theorem}

\begin{figure}[t]
\centering
\begin{minipage}{0.49\textwidth}\centering
\includegraphics[width=\linewidth]{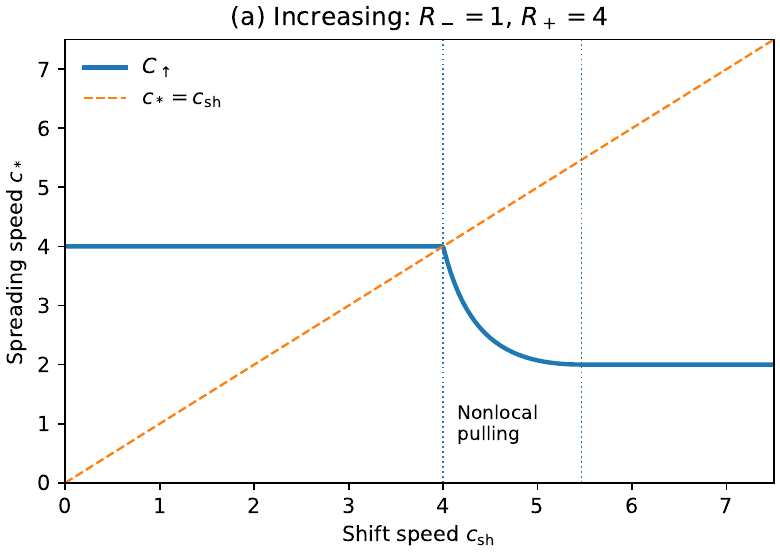}
\end{minipage}\hfill
\begin{minipage}{0.49\textwidth}\centering
\includegraphics[width=\linewidth]{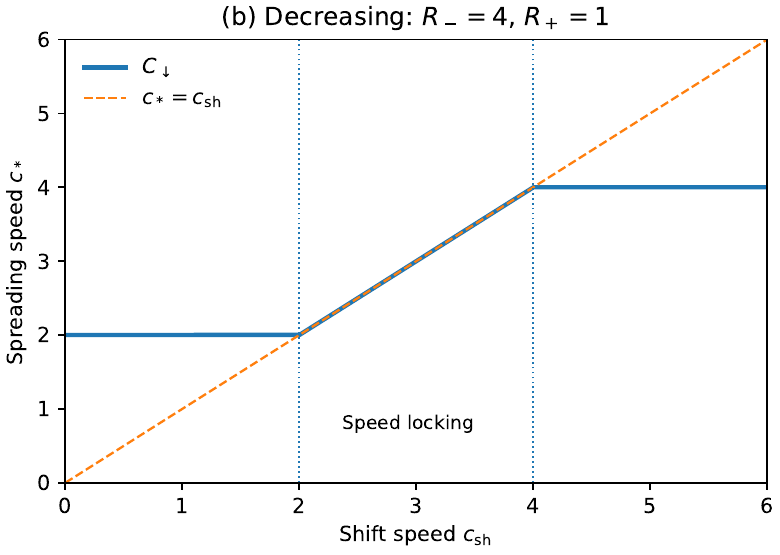}
\end{minipage}
\caption{The exact spreading  speed formulas. Left: $R_-=1$, $R_+=4$; the
intermediate zone is nonlocally pulled. Right: $R_-=4$, $R_+=1$;
the intermediate zone is speed locking. The solid blue  curves describe the spreading speed formula, while,  the dashed orange  curves indicate
$\cstar=\csh$.}
\label{fig:speed-curves}
\end{figure}

\begin{remark}
In this  case, $\cstar=\csh$ when $2\sqrt{R_+}<\csh<2\sqrt{R_-}$.  This equality does not imply that a level set remains at a bounded distance from the interface $A(t)$ with Assumption~\ref{ASS-A} alone. 
 Figure~\ref{fig:speed-curves} illustrates the dependence of the spreading speeds in Theorems~\ref{Thm-inc} and \ref{Thm-dec} on $\csh$.  
\end{remark}

\begin{remark}
If  $r(t,z)$ is between $r_-(t)$ and $r_+(t)$ in Assumption~\ref{ASS-limits} and  $R_-=R_+=R>0$, then $r(t,z)\equiv  r_-(t) = r_+(t)$. The spreading speed is $c_*= 2\sqrt{R}$.
\end{remark}

\section{The stepwise problem}\label{Sec-step}
We first study an auxiliary problem with a piecewise $T$-periodic step reaction. Throughout this section, $f_-$ and $f_+$ denote continuous $T$-periodic functions on $\R\times[0,M]$ satisfying common boundedness and Lipschitz estimates together with $r_\pm(t):=f_\pm(t,0)$ and 
$$
 r_\pm(t)-Ku\leq f_\pm(t,u)\leq r_\pm(t),\quad
 \quad f_\pm(t,M)\leq-\sigma_0,\qquad  R_\pm=\av{r_\pm}>0. 
$$
Assume the shift is linear,  $A(t)=c_1t$ with $c_1>0$.  Define
\begin{equation}\label{stepg}
g (t,z,u):=
 \begin{cases}
 f_-(t,u),&z<0,\\
 f_+(t,u),&z\geq0.
 \end{cases}
\end{equation}
We consider
\begin{equation}\label{eq-step}
 \begin{cases}
 u_t=u_{xx}+u g(t,x-c_1t,u),&t>0,\ x\in\R,\\
 u(0,x)=u_0(x),&x\in\R.
 \end{cases}
\end{equation}
Note  that \eqref{eq-step} has a unique  global bounded mild solution and maximum principle holds.   Since the line $x=c_1 t$ has zero measure, changing the value of  $g$ at $x=c_1 t$ does not change the mild solution.

For $0<\alpha<\beta$ and $c>0$, define
\begin{equation}\label{eq:C-inc}
 C_\uparrow(c;\alpha,\beta):=
 \begin{cases}
  2\sqrt\beta,
  \quad   &0<c\leq2\sqrt\beta,\\[1mm]
  \displaystyle \frac{c}{2}-\sqrt{\beta-\alpha}
  +\frac{\alpha}{\frac{c}{2}-\sqrt{\beta-\alpha}},
    &2\sqrt\beta<c<2\bigl(\sqrt{\beta-\alpha}+\sqrt\alpha\bigr),\\[3mm]
  2\sqrt\alpha,
    &c\geq2\bigl(\sqrt{\beta-\alpha}+\sqrt\alpha\bigr).
 \end{cases}
\end{equation}
For $a>b>0$, define
\begin{equation}\label{eq:C-dec}
 C_\downarrow(c;a,b):=
 \begin{cases}
  2\sqrt b,\quad &0<c\leq2\sqrt b,\\
  c,&2\sqrt b<c<2\sqrt a,\\
  2\sqrt a,&c\geq2\sqrt a.
 \end{cases}
\end{equation}

The maps $(c,\alpha,\beta)\mapsto C_\uparrow(c;\alpha,\beta)$ and
$(c,a,b)\mapsto C_\downarrow(c;a,b)$ are jointly continuous on
$\{c>0,\ 0<\alpha<\beta\}$ and $\{c>0,\ a>b>0\}$, respectively.
Indeed, for \eqref{eq:C-inc}, at $c=2\sqrt\beta$,
$$
 \lambda(c):=\frac c2-\sqrt{\beta-\alpha}
 =\sqrt\beta-\sqrt{\beta-\alpha},
 \qquad
 \lambda + \frac{\alpha}{\lambda} 
 =2\sqrt\beta.
$$
   At
$c=2(\sqrt{\beta-\alpha}+\sqrt\alpha)$, one has
$\lambda=\sqrt\alpha$, and the middle expression equals
$2\sqrt\alpha$.  The corresponding boundary identities for
\eqref{eq:C-dec} are immediate.  In the middle interval, since $\lambda(c) < \sqrt{\alpha}$, we have 
$$
 \frac{d}{dc} \left(\lambda(c)+\frac{\alpha}{\lambda(c)}\right)
 =\frac12\left(1-\frac{\alpha}{\lambda(c)^2}\right)<0.
$$
Hence $C_\uparrow$ is nonincreasing in $c$, whereas $C_\downarrow$ is
nondecreasing in $c$.

\begin{theorem}\label{Thm-step}
Let $c_1>0$ and let $u_0\in C (\R,[0,M])$ be nonzero and compactly supported.  Then:
\begin{itemize}
 \item[(i)] If $f_-(t,u)\leq f_+(t,u)$ for every
 $(t,u)\in[0,T]\times[0,M]$ and
 $0<R_-<R_+$, the spreading speed of \eqref{eq-step} is
$$
 C_\uparrow(c_1;R_-,R_+).
$$
 \item[(ii)] If $f_-(t,u)\geq f_+(t,u)$ for every
 $(t,u)\in[0,T]\times[0,M]$ and
 $R_->R_+>0$, the spreading speed of \eqref{eq-step} is
$$
 C_\downarrow(c_1;R_-,R_+).
$$
\end{itemize}
 The conclusion still holds if the jump line $x=c_1t$ is replaced by
 $x=c_1t+D$ for any constant $D\in \R$.
\end{theorem}

The last conclusion can be  obtained by the spatial shifting. Here  we only prove the case of $x=c_1 t$.
To prove Theorem~\ref{Thm-step}, we first prove three lemmas which are elementary blocks  in constructing sub and supersolutions.  Let $b(t)$ be a continuous $T$-periodic function with $\av{b}>0$, and set
$$
 \theta_b(t):=\exp\left\{\int_0^t\bigl(b(s)-\av b\bigr)\,ds\right\}.
$$
Note that $\theta_b$ is positive and $T$-periodic.

\begin{lemma}\label{lem:periodic-plateau}
Assume that
$$
 F(t,x,u)\geq b(t)-Ku,
 \qquad (t,x)\in\Omega,\quad 0\leq u\leq M,
$$
where $\Omega\subset[0,\infty)\times\R$.  If $\eta>0$ satisfies
$$
 0<\eta\max_{t\in[0,T]}\theta_b(t)
 \leq\min\left\{M,\frac{\av b}{2K}\right\},
$$
then $v(t,x):=\eta   \theta_b(t)$ satisfies
$$
 v_t-v_{xx}-vF(t,x,v)\leq0,\qquad (t,x)\in\Omega.
$$
\end{lemma}

\begin{proof}
Since $\theta_b'=(b-\av b)\theta_b$, we have 
\begin{equation*}
 \begin{split}
 v_t-v_{xx}-vF(t,x,v)
 &\leq (b(t)-\av b)\eta\theta_b(t)
       -b(t)\eta\theta_b(t)+K\eta^2\theta_b^2(t)\\
 &=\eta\theta_b(t)\bigl(-\av b+K\eta\theta_b(t)\bigr)\\
 &\leq-\frac{\eta\av b\,\theta_b(t)}2<0.
 \end{split}
\end{equation*}
The proof is complete.
\end{proof}

\begin{lemma}\label{Lem-super}
For $\mu>0$, define
$$
\phi_{b,\mu}(t,x):=\theta_b(t)
\exp\left\{-\mu x+(\mu^2+\av b)t\right\}.
$$
Then
$$
 (\phi_{b,\mu})_t-(\phi_{b,\mu})_{xx}
 -b(t)\phi_{b,\mu}=0\qquad\text{in }\R^2.
$$
\end{lemma}
\begin{proof}
Direct computation gives
$$
 (\phi_{b,\mu})_t=(b(t)+\mu^2)\phi_{b,\mu},
 \qquad
 (\phi_{b,\mu})_{xx}=\mu^2\phi_{b,\mu}.
$$
The proof is complete
\end{proof}

For $\alpha\geq0$, $c\geq0$, $D\in\R$ and $l>0$, define
\begin{equation}\label{eq-sub}
 \psi_{b,\alpha,c,l,D}(t,x):=
 \begin{cases}
 \theta_b(t)\exp\left\{-\alpha t-\frac{c}{2}(x-ct-D)\right\}
 \sin\left(\frac{x-ct-D}{l}\right),
 &0<\dfrac{x-ct-D}{l}<\pi,\\
 0,&\text{otherwise}.
 \end{cases}
\end{equation}

\begin{lemma}\label{Lem-sub}
Let $F=F(t,x, u)$ be measurable in $(t,x)$ for each $u \in [0, M]$, bounded by $B$, uniformly Lipschitz continuous in $s$ with constant $L$, and $F(t,x,M )\leq 0$. 
Assume that, on the support of
$\psi=\psi_{b,\alpha,c,l,D}$,
$$
 F(t,x,u)\geq b(t)-Ku,\qquad 0\leq u\leq M.
$$
If
$$
 \frac{c^2}{4}+\frac1{l^2}<\av b+\alpha,
$$
and $\eta>0$ is sufficiently small, 
then $0\leq\eta\psi\leq M$ and $\eta\psi$ satisfies, for $t>0$ and
$0<\frac{x-ct-D}{l}<\pi$,
\begin{equation}\label{subeq}
 (\eta\psi)_t-(\eta\psi)_{xx}
 -\eta\psi F(t,x,\eta\psi)\leq0.
\end{equation}
The function $\eta\psi $  is a
generalized subsolution.
\end{lemma}

\begin{proof}
For $t>0$ and $0<\frac{x-ct-D}{l}<\pi$, we have 
\begin{equation*}
\psi_t
=\left(b(t)-\av b-\alpha+\frac{c^2}{2}\right)\psi
-\frac{c}{l}\theta_b(t)e^{-\alpha t-\frac c2(x-ct-D)}
\cos\left(\frac{x-ct-D}{l}\right),
\end{equation*}
\begin{equation*}
\psi_{xx}
=\left(\frac{c^2}{4}-\frac1{l^2}\right)\psi
-\frac{c}{l}\theta_b(t)e^{-\alpha t-\frac c2(x-ct-D)}
\cos\left(\frac{x-ct-D}{l}\right),
\end{equation*}
and 
\begin{equation}\label{subproof}
\psi_t-\psi_{xx}-b(t)\psi
=\left(\frac{c^2}{4}+\frac1{l^2}-\av b-\alpha\right)\psi.
\end{equation}
Set
$$
 \delta:=-\frac{c^2}{4}-\frac1{l^2}+\av b+\alpha>0.
$$
Since $F(t,x,u)\geq b(t)-Ku$ and \eqref{subproof}, we  obtain, for small $\eta>0$, 
$$
 (\eta\psi)_t-(\eta\psi)_{xx}
 -\eta\psi F(t,x,\eta\psi)
 \leq\eta\psi(-\delta+K\eta\psi).
$$
Choose $\eta$ so small that
$$
 \eta\|\psi\|_{L^\infty([0,\infty)\times\R)}\leq M,
 \qquad
 K\eta\|\psi\|_{L^\infty([0,\infty)\times\R)}\leq\delta.
$$
Then $\eta\psi$ satisfies \eqref{subeq}.  
Hence $\eta \psi$  is a generalized
subsolution.  The proof is complete.
\end{proof}

We shall also use the following lower estimate.

\begin{lemma}\label{Lem-lowerbound}
Let $F$ satisfy the hypotheses of Lemma \ref{Lem-sub}, and let $v$ be the solution starting from nontrivial compactly supported initial data, with $0\leq v\leq M$, of
$$
 v_t=v_{xx}+vF(t,x,v).
$$
Suppose that there is a continuous $T$-periodic function $b$ with
$\av{b}>0$ such that
$$
 F(t,x,s)\geq b(t)-Ks,
 \qquad t\geq0,\quad x\in\R,\quad0\leq s\leq M.
$$
Then
$$
 \liminf_{t\to\infty}\inf_{0\leq x\leq qt}v(t,x)>0
 \qquad  \forall  \,  0<q<2\sqrt{\av{b}}.
$$
The left endpoint $0$ may be replaced by any fixed   $a\in\R$.
\end{lemma}

\begin{proof}
By maximum principle, we know that $v(t,x)>0$ for all $t>0$ and $x\in \R$. 
Fix
$$
 0<q_0<q<2\sqrt{\av{b}}
$$
and choose $l_0,l_q>0$ such that
$$
 l_0^{-2}<\av{b},
 \qquad
 \frac{q^2}{4}+l_q^{-2}<\av{b}.
$$
Let $0<\eta<1$.  In \eqref{eq-sub}, take $\alpha=0$ and set
$$
 D_0:=-1-l_0(\pi-\arcsin\eta),
$$
$$
 \Psi_0(t,x):=\eta\psi_{b,0,0,l_0,D_0}(t,x),
 \qquad P(t):=\eta^2\theta_b(t).
$$
Then $\Psi_0(t,-1)=P(t)$ and $(\Psi_0)_x(t,-1)<0$.  Next set
$$
 h_q(y):=e^{-qy/2}\sin(y/l_q),
 \qquad 0\leq y\leq\pi l_q.
$$
For sufficiently small $\eta$, the equation $h_q(y)=\eta$ has a first
root $y_q$ on the increasing part of $h_q$.  Set
$$
 \Psi_q(t,x):=\eta\psi_{b,0,q,l_q,0}(t,x),
 \qquad \xi_q(t):=qt+y_q.
$$
Thus $\Psi_q(t,\xi_q(t))=P(t)$ and
$(\Psi_q)_x(t,\xi_q(t))>0$.

Define
\begin{equation}\label{eq:lowerbound-subsolution}
 \underline v(t,x)=
 \begin{cases}
 0,&x\leq D_0,\\
 \Psi_0(t,x),&D_0<x<-1,\\
 P(t),&-1\leq x\leq\xi_q(t),\\
 \Psi_q(t,x),&\xi_q(t)<x<qt+\pi l_q,\\
 0,&x\geq qt+\pi l_q.
 \end{cases}
\end{equation}
Choose $\eta$ sufficiently small such  that
$0\leq\underline v\leq M$.  The functions $\Psi_0$, $P$ and
$\Psi_q$ then satisfy the subsolution inequality by Lemmas
\ref{lem:periodic-plateau} and \ref{Lem-sub}.  Note that 
\begin{align*}
&\underline{v}_x(t, D_0+)  -\underline{v}_x(t, D_0-) = \frac{\eta\theta_b(t)}{l_0} \geq 0,  \\
&\underline{v}_x(t, (-1)+)  -\underline{v}_x(t, (-1)-) =\frac{\eta\theta_b(t)}{l_0}\sqrt{1-\eta^2}\ge 0,\\
&\underline{v}_x(t, \xi_q(t)+)  -\underline{v}_x(t, \xi_q(t)-) =  \eta\theta_b(t)h_q'(y_q)\geq 0, \\
&\underline{v}_x(t, (qt+ \pi l_q)+)  -\underline{v}_x(t, (qt+ \pi l_q)-) = \frac{\eta\theta_b(t)}{l_q}e^{-q\pi l_q/2}\ge 0. 
\end{align*}
  Hence
\eqref{eq:lowerbound-subsolution} is a generalized subsolution.

Fix $t_0>0$.   We can choose $\eta>0$ sufficiently small such that $\underline{v}(t_0, \cdot) \leq v(t_0, \cdot)$.
 Since
$q_0t<\xi_q(t)$ for all sufficiently large $t$, comparison principle  yields
$$
 \liminf_{t\to\infty}\inf_{0\leq x\leq q_0t}v(t,x)
 \geq\eta^2\min_{[0,T]}\theta_b>0.
$$
Since $q_0<2\sqrt{\av b}$ is arbitrary, we obtain  the conclusion.

For a fixed $a\in\R$, replace the  function $\Psi_0(t,x)$ by
$$
 \Psi_0^a(t,x):=\Psi_0(t,x-a),
 \qquad D_0^a:=D_0+a,
$$
and replace the matching point $-1$ by $a-1$.   We can show 
$$
 \liminf_{t\to\infty}\inf_{a \leq x\leq qt}v(t,x)>0
 \qquad  \forall  \,  0< q<2\sqrt{\av{b}}.
$$
The proof is complete.
\end{proof}

We now prove Theorem \ref{Thm-step}.

\subsection{Increasing stepwise profile}
\begin{proof} [Proof of Theorem \ref{Thm-step}\textit{(i)}]
We first prove the upper estimates.  Since $g(t,z,s)\leq f_+(t,s)\leq
r_+(t)$, the function
$$
 \overline u(t,x)=\min\{M,B\phi_{r_+,\sqrt{R_+}}(t,x)\}
$$
is a generalized supersolution.  Along the interface where
$B\phi_{r_+,\sqrt{R_+}}=M$, the spatial derivative gap equals
$-\sqrt{R_+}M$.  After increasing $B$,  we have $\overline{u}(0,\cdot) \geq u_0$,
and for every $c>2\sqrt{R_+}$,
$$
 \sup_{x\geq ct}\overline u(t,x)
 \leq C e^{-\sqrt{R_+}(c-2\sqrt{R_+})t}\longrightarrow0,  \quad \text{ as } t\to +\infty.
$$
Hence $\overline c(u)\leq2\sqrt{R_+}$.

For the other two ranges, take $0<\lambda<\mu$ and put
$$
 U_-(t,x):=B_-\phi_{r_-,\lambda}(t,x),
 \qquad
 U_+(t,x):=B_+\phi_{r_+,\mu}(t,x).
$$
Their intersection curve is
\begin{equation}\label{eq:inc-equality-curve}
 X(t)=c_{\lambda,\mu}t+\frac{1}{\mu-\lambda}
 \log\frac{B_+\theta_{r_+}(t)}{B_-\theta_{r_-}(t)},
 \qquad
 c_{\lambda,\mu}:=\lambda+\mu+
 \frac{R_+-R_-}{\mu-\lambda}.
\end{equation}

When $c_{\lambda,\mu}=c_1$, we choose  $B_+/B_-$ such  that $X(t)<c_1t$ for all $t\geq0$, and then multiply  $B_-$ and $B_+$ by the same large constant  so that
\begin{equation}\label{eq:inc-general-super}
 \overline u(t,x):=\min\{M,U_-(t,x),U_+(t,x)\}
\end{equation}
satisfies $\overline{u}(0, \cdot) \geq u_0(\cdot)$.   Then we can  calculate the negative derivative gap. Hence \eqref{eq:inc-general-super} is a generalized supersolution.

Assume 
$$
 2\sqrt{R_+}<c_1<
 2\bigl(\sqrt{R_+-R_-}+\sqrt{R_-}\bigr).
$$
 Set
 \begin{equation}\label{eq:step-lambda-star}
  \lambda_*:=\frac{c_1}{2}-\sqrt{R_+-R_-},
  \qquad
  c_{\rm nlp}:=\lambda_*+\frac{R_-}{\lambda_*}.
 \end{equation}
 Taking $\lambda=\lambda_*$ and $\mu=c_1/2$ in
 \eqref{eq:inc-equality-curve} gives $c_{\lambda,\mu}=c_1$.  Moreover,
 $$
  \lambda_*c_1-\lambda_*^2-R_-
  =\frac{c_1^2}{4}-R_+>0.
 $$
 
 If $c>c_{\rm nlp}$, then
 $$
 \lim_{t\to +\infty} \sup_{x\geq ct}  \overline u(t,x)\leq  \lim_{t\to +\infty} \sup_{x\geq ct}  U_-(t,x)
  \leq   \lim_{t\to +\infty}  C\exp\{-\lambda_*(c-c_{\rm nlp})t\} = 0.
 $$
 Hence
 $$
  \overline c(u)\leq c_{\rm nlp}.
 $$
 
 If  $
  c_1\geq2\bigl(\sqrt{R_+-R_-}+\sqrt{R_-}\bigr)$, we 
take  $\lambda=\sqrt{R_-}$ and    choose $h>0$ such that
 $$
  h+\frac{R_+-R_-}{h}=c_1-2\sqrt{R_-},
  \qquad \mu=\sqrt{R_-}+h.
 $$
The existence of $h$ is  because
 $c_1-2\sqrt{R_-}\geq2\sqrt{R_+-R_-}$.
Then  $c_{\lambda,\mu}=c_1$ in \eqref{eq:inc-equality-curve}.  The supersolution
 \eqref{eq:inc-general-super} satisfies, for every $c>2\sqrt{R_-}$,
 $$
  \lim_{t\to +\infty}  \sup_{x\geq ct}\overline u(t,x)
  \leq C  \lim_{t\to +\infty}  e^{-\sqrt{R_-}(c-2\sqrt{R_-})t} =0.
 $$
Thus,  $
  \overline c(u)\leq2\sqrt{R_-}.$
 We obtain 
 \begin{equation}\label{eq:inc-upper-bound}
  \overline c(u)\leq C_\uparrow(c_1;R_-,R_+).
 \end{equation}

 We turn to the lower estimates.  The increasing order gives
 $$
  g(t,z,s)\geq f_-(t,s)\geq r_-(t)-Ks
 $$
 for every $z\in\R$.  It follows from  Lemma \ref{Lem-lowerbound} that
 $$
  \underline c(u)\geq2\sqrt{R_-}.
 $$
 
 Assume  $c_1<2\sqrt{R_+}$.  Fix
 $0<q_0<2\sqrt{R_+}$ and  let 
 $$
  \max\{q_0,c_1\}<q<2\sqrt{R_+}.
 $$
 We  choose $l_0,l_c,l_q>0$ such that
 \begin{equation}\label{eq:inc-low-radii}
  l_0^{-2}<R_-,
  \qquad \frac{c_1^2}{4}+l_c^{-2}<R_+,
  \qquad \frac{q^2}{4}+l_q^{-2}<R_+.
 \end{equation}
 For $0<\eta<1$, set
 $$
  P_-(t)=\eta^2\theta_{r_-}(t),
  \qquad P_+(t)=\eta^2\theta_{r_+}(t).
 $$
Let
 $$
  D_0:=-1-l_0(\pi-\arcsin\eta),
 $$
 such  that
 $$
  \Psi_0(t,x):=\eta\psi_{r_-,0,0,l_0,D_0}(t,x)
 $$
equals to  $P_-$ at $x=-1$ on the decreasing part of $\Psi_0$.  Define 
 $$
  \Psi_c(t,x):=\eta\psi_{r_+,0,c_1,l_c,1}(t,x).
 $$
 For $0\leq y\leq\pi l_c$, define
 $$
  h_c(y):=e^{-c_1y/2}\sin(y/l_c),
 \qquad m_c:=\max_{[0,\pi l_c]}h_c>0.
 $$
 There is a unique $\bar y_c\in(0,\pi l_c)$ 
 such that $h_c(\bar{y}_c) =m_c$.  
 The function $h_c(y)$  is   increasing on $[0,\bar y_c]$ and
 decreasing on $[\bar y_c,\pi l_c]$.  Since
 $\theta_{r_-}/\theta_{r_+}$ is positive and periodic, it is bounded
 above and below by positive constants.  
 There exists
 $\eta_0\in(0,1)$ such that, for every $0<\eta\leq\eta_0$,
 $$
  \eta\frac{\theta_{r_-}(t)}{\theta_{r_+}(t)}<m_c,
  \qquad \eta<m_c,
  \qquad  \forall t\in\R,
 $$
 and 
 $$\eta<\max_{[0,\pi l_q]}e^{-qy/2}\sin(y/l_q). $$
 For such $\eta$, let $y_c^-(t)\in(0,\bar y_c)$ and
 $y_c^+\in(\bar y_c,\pi l_c)$ be the unique solutions of
 $$
  h_c(y_c^-(t))
  =\eta\frac{\theta_{r_-}(t)}{\theta_{r_+}(t)},
  \qquad h_c(y_c^+)=\eta.
 $$
 Set
 $$
  \xi_c^-(t):=c_1t+1+y_c^-(t),
  \qquad
  \xi_c^+(t):=c_1t+1+y_c^+.
 $$
Note  that $y_c^-$ is $C^1$ and
 $T$-periodic.  The definitions yield  
 $$
  \Psi_c(t,\xi_c^-(t))=P_-(t),\qquad
  \Psi_c(t,\xi_c^+(t))=P_+(t),
 $$
 and
 \begin{equation}\label{eq:inc-low-cap-derivatives}
  (\Psi_c)_x(t,\xi_c^-(t))>0,\qquad
  (\Psi_c)_x(t,\xi_c^+(t))<0.
 \end{equation}
 Both $\xi_c^-(t)-c_1t$ and $\xi_c^+(t)-c_1t$ are bounded
 $T$-periodic functions.
 Finally, let
 $$
  \Psi_q(t,x):=\eta\psi_{r_+,0,q,l_q,0}(t,x).
 $$
 Writing $h_q(y)=e^{-qy/2}\sin(y/l_q)$, let $y_q$ be the unique
 solution of $h_q(y_q)=\eta$ on the increasing part  of $h_q$, and
 set $\xi_q(t):=qt+y_q$.  Then
 $$
  \Psi_q(t,\xi_q(t))=P_+(t),\qquad
  (\Psi_q)_x(t,\xi_q(t))>0.
 $$
 Choose $t_0>0$ so large that
 $$
  c_1t+1+\pi l_c<qt,\qquad t\geq t_0.
 $$
 Since
 $$
  \xi_c^+(t)<c_1t+1+\pi l_c,
  \qquad \xi_q(t)>qt,
 $$
 we have, for every $0<\eta\leq\eta_0$,
 $$
  -1<\xi_c^-(t)<\xi_c^+(t)<\xi_q(t),
  \qquad t\geq t_0.
 $$
 Fix $\eta\in(0,\eta_0]$.
 For $t\geq t_0$, define
 \begin{equation}\label{eq:inc-low-subsolution}
  \underline u(t,x)=
  \begin{cases}
  0,&x\leq D_0,\\
  \Psi_0(t,x),&D_0<x<-1,\\
  P_-(t),&-1\leq x\leq\xi_c^-(t),\\
  \Psi_c(t,x),&\xi_c^-(t)<x<\xi_c^+(t),\\
  P_+(t),&\xi_c^+(t)\leq x\leq\xi_q(t),\\
  \Psi_q(t,x),&\xi_q(t)<x<qt+\pi l_q,\\
  0,&x\geq qt+\pi l_q.
  \end{cases}
 \end{equation}
 
 Define operator 
 $$
  \mathcal L_g[w]:=w_t-w_{xx}
  -w g(t,x-c_1t,w).
 $$
 The increasing property gives $g(t,z,w)\geq r_-(t)-Kw$ for  $z\in \R$.  On $x>c_1t$, the  assumption of $g$ implies 
 $g(t,x-c_1t, w)=f_+(t,w)\geq r_+(t)-Kw$.  Set
 $$
  \gamma_0:= R_{-} -l_0^{-2},\quad
  \gamma_c:=R_+-\frac{c_1^2}{4}-l_c^{-2},\quad
  \gamma_q:=R_+-\frac{q^2}{4}-l_q^{-2}.
 $$
 All three constants are positive by \eqref{eq:inc-low-radii}.  On the
 corresponding domain,
 \begin{align*}
  \mathcal L_g[\Psi_0]&\leq-\gamma_0\Psi_0+K\Psi_0^2,&
  \mathcal L_g[P_-]&\leq-R_-P_-+KP_-^2,\\
  \mathcal L_g[\Psi_c]&\leq-\gamma_c\Psi_c+K\Psi_c^2,&
  \mathcal L_g[P_+]&\leq-R_+P_++KP_+^2,\\
  \mathcal L_g[\Psi_q]&\leq-\gamma_q\Psi_q+K\Psi_q^2.
 \end{align*}
 By direct computation, we  also know that the spatial derivative gap is  positive at each intersection line.   Hence,   \eqref{eq:inc-low-subsolution} is a generalized subsolution.

 Fix $t_0$.  We can choose $\eta>0$ sufficiently small such that $\underline{u}(t_0, \cdot) \leq u(t_0, \cdot)$. 
 Note that 
 $$
   \Psi_c(t,x)\geq\min\{P_-(t),P_+(t)\}
   \quad\text{for }\xi_c^-(t)\leq x\leq\xi_c^+(t),
 $$
 and $q_0t<\xi_q(t)$ for all large $t$.  Hence the subsolution $\underline u$ has a positive lower
  bound on $0\leq x\leq q_0t$.
 By comparison principle, we obtain 
 $$
  \liminf_{t\to\infty}\inf_{0\leq x\leq q_0t}u(t,x)
  \geq\eta^2
  \min\left\{\min_{[0,T]}\theta_{r_-},
              \min_{[0,T]}\theta_{r_+}\right\}>0.
 $$
 Thus 
 $$
  \underline c(u)\geq2\sqrt{R_+} \quad \text{ when } c_1 <2\sqrt{R_+}.
  $$

  It remains to prove the lower estimate in the intermediate range,  $2\sqrt{R_+} \le c_1 < 2(\sqrt{R_+ -R_{-}} + \sqrt{R_-})$.  Let 
  $\lambda_*$ and $c_{\rm nlp}$ be defined as in \eqref{eq:step-lambda-star}.  Fix $0<c_0<c_{\rm nlp}$.  In this range,
  $0<\lambda_*<\sqrt{R_-}$.  Choose
  $$
   0<\delta<
   \min\left\{R_-,\lambda_*(c_{\rm nlp}-c_0), {R_--\lambda_*^2}\right\}.
  $$
  Then
  $$
   \omega_\delta:=\lambda_*+\frac{R_--\delta}{\lambda_*}
   =c_{\rm nlp}-\frac{\delta}{\lambda_*}
   >\max\{c_0,2\lambda_*\}.
  $$
 Set
  $$
   d_\delta:=\lambda_*(c_1-\omega_\delta)
   =\frac{c_1^2}{4}-R_++\delta>0.
  $$
Let $0<\nu<\omega_\delta-2\lambda_*$ and 
  $\xi=x-\omega_\delta t$.  Define
  $$
   \rho(t,x):=\theta_{r_-}(t)
   \left[e^{-\lambda_*\xi}-e^{-(\lambda_*+\nu)\xi}\right]_+.
  $$
  In $\xi>0$,
  \begin{equation}\label{eq:rho-jin}
   \rho_t-\rho_{xx}-r_-(t)\rho
   =-\delta\rho
   -\nu(\omega_\delta-2\lambda_*-\nu)\theta_{r_-}(t)
   e^{-(\lambda_*+\nu)\xi}
   \leq -\delta\rho.
  \end{equation}
  Choose $l>0$ so that $l^{-2}<\delta/2$ and set
  $$
   \widehat\psi(t,x):=\psi_{r_+,d_\delta,c_1,l,1}(t,x).
  $$
  The support of $\widehat\psi$ lies in $x>c_1t$, and
  \begin{equation}\label{eq:psi-hat-jin}
   \widehat\psi_t-\widehat\psi_{xx}-r_+(t)\widehat\psi
   =(-\delta+l^{-2})\widehat\psi.
  \end{equation}
  
Let $l_0$ satisfy  $l_0^{-2}<R_-$.  For $0<\eta<1$, denote 
  $$
   D_0=-1-l_0(\pi-\arcsin\eta),
  $$
  $$
   w_{-1}(t,x)=\eta^2\psi_{r_-,0,0,l_0,D_0}(t,x),
   \quad w_0(t,x)=\eta^3\theta_{r_-}(t),
   \quad w_1(t,x)=\eta^2\rho(t,x),
   \quad w_2(t,x)=\eta\widehat\psi(t,x).
  $$
  Then $w_{-1}(t, -1)=w_0(t, -1)$.  Let
  $$
   p(\xi):=e^{-\lambda_*\xi}-e^{-(\lambda_*+\nu)\xi},
   \qquad \xi\geq0.
  $$
  The function $p$ attains its  maximum at
  $$
   \xi_p:=\frac1\nu\log\frac{\lambda_*+\nu}{\lambda_*}>0.
  $$
  There exists $\eta_*\in(0,1)$ such that, for every
  $0<\eta\leq\eta_*$, the equation $p(\xi_1)=\eta$ has a unique solution
  $\xi_1=\xi_1(\eta)\in(0,\xi_p)$ on the increasing part.  Define
  $$
   X_1(t)=\omega_\delta t+\xi_1
  $$
  so that $w_0=w_1$ at $x=X_1(t)$  and
  $$
   (w_1)_x(t,X_1(t))-(w_0)_x(t,X_1(t))
   =\eta^2\theta_{r_-}(t)p'(\xi_1)>0.
  $$
  
  Choose $y_*\in(0,\pi l)$ so that
  $e^{-c_1y/2}\sin(y/l)$ is strictly increasing on $[0,y_*]$. Let 
  $t_0>0$  be sufficiently large such that, for all $t\ge t_0$, 
  \begin{equation}\label{eq:pulling-uniform-starting-time}
    \omega_\delta t+\xi_p<c_1t+1, \quad 
    e^{-\nu((c_1-\omega_\delta)t+1)}<\frac{\lambda_*}{\lambda_*+\nu}.
  \end{equation}
  At
  $x=c_1t+1+y$,
  \begin{align*}
   w_1(t,x)
   &=\eta^2\theta_{r_-}(t)e^{-d_\delta t}e^{-\lambda_* (1+y)}
   \left[1-e^{-\nu((c_1-\omega_\delta)t+1+y)}\right],\\
   w_2(t,x)
   &=\eta\theta_{r_+}(t)e^{-d_\delta t}
   e^{-c_1y/2}\sin(y/l).
  \end{align*}
  The second inequality in
  \eqref{eq:pulling-uniform-starting-time} implies that the map
  $y\mapsto w_1(t,c_1t+1+y)$ is  decreasing on $[0,y_*]$ for
  $t\geq t_0$, whereas $y\mapsto w_2(t,c_1t+1+y)$ is 
  increasing there.  
  At $y=0$,
  $$
   w_1(t,c_1t+1)>0=w_2(t,c_1t+1).
  $$
  Moreover, periodicity and positivity
  of $\theta_{r_\pm}$ give a constant $C_0>0$, independent of $t$ and
  $\eta$, such that
  $$
   \frac{w_1(t,c_1t+1+y_*)}{w_2(t,c_1t+1+y_*)}
   \leq C_0\eta.
  $$
Decreasing $\eta_*$ such  that $C_0\eta_*<1$, we obtain a unique
  $y_2(t)=y_2(t;\eta)\in(0,y_*)$ for which
  $$
   X_2(t)=c_1t+1+y_2(t),
   \qquad w_1(t,X_2(t))=w_2(t,X_2(t)),
  $$
  and $(w_2)_x(t,X_2(t))-(w_1)_x(t,X_2(t))>0$. By the implicit function theorem, $y_2$ is $C^1$ for $t\ge t_0$. Moreover, $y_2$ is bounded in $(0,y_*)$ for $t\geq t_0$.
  Let $X_3(t)=c_1t+1+\pi l$.  By
  \eqref{eq:pulling-uniform-starting-time}, for $t\geq t_0$,
  $$
   D_0<-1<0<X_1(t)<X_2(t)<X_3(t).
  $$
  Fix $\eta\in(0,\eta_*]$ small enough.  
  Define
  \begin{equation}\label{eq:pulling-subsolution-jin}
   \underline u(t,x)=
   \begin{cases}
   0,&x\leq D_0,\\
   w_{-1}(t,x),&D_0<x<-1,\\
   w_0(t),&-1\leq x\leq X_1(t),\\
   w_1(t,x),&X_1(t)<x\leq X_2(t),\\
   w_2(t,x),&X_2(t)<x<X_3(t),\\
   0,&x\geq X_3(t).
   \end{cases}
  \end{equation}
  Set
  $$
   \gamma_{-1}:=R_- -l_0^{-2}>0,\qquad
   \gamma_2:=\delta-l^{-2}>\frac{\delta}{2}.
  $$
  The increasing order, \eqref{eq:rho-jin}, and
  \eqref{eq:psi-hat-jin} give the following estimates on each domain:
  \begin{align*}
   \mathcal L_g[w_{-1}]
   &\leq-\gamma_{-1}w_{-1}+Kw_{-1}^2,\\
   \mathcal L_g[w_0]
   &\leq-R_-w_0+Kw_0^2,\\
   \mathcal L_g[w_1]
   &\leq-\delta w_1+Kw_1^2,\\
   \mathcal L_g[w_2]
   &\leq-\gamma_2w_2+Kw_2^2.
  \end{align*}
 The derivative gap is positive at each intersection lines.  Thus \eqref{eq:pulling-subsolution-jin} is a generalized subsolution.
  

  Let $t_0$ be fixed in   \eqref{eq:pulling-uniform-starting-time}. We can choose $\eta $ small enough such that 
  $\underline u(t_0,\cdot)\leq u(t_0,\cdot)$.  
The comparison principle  yields 
  $$
   \liminf_{t\to\infty}\inf_{0\leq x\leq c_0t}u(t,x)
   \geq\eta^3\min_{[0,T]}\theta_{r_-}>0,
  $$
 By the arbitrariness of $c_0 \in (0, c_{\rm nlp})$, we have 
  $$
   \underline c(u)\geq c_{\rm nlp}.
  $$
 
  Combining these lower bounds with \eqref{eq:inc-upper-bound}, we  complete the proof of 
  Theorem \ref{Thm-step}\textit{(i)}.

\end{proof}

\subsection{Decreasing stepwise profile}
\begin{proof}[Proof of Theorem \ref{Thm-step} \textit{(ii)}]
We first prove the upper estimates.  Suppose that
$0<c_1\leq2\sqrt{R_+}$.  Choose $B>0$  large enough such that
$$
 \xi_1(t):=2\sqrt{R_+}\,t+
 \frac1{\sqrt{R_+}}\log\frac{B\theta_{r_+}(t)}M>c_1t,
 \qquad  \forall  t\geq0,
$$
and
$$
 \mathcal U_1(t,x):=
 \min\{M,B\phi_{r_+,\sqrt{R_+}}(t,x)\} \quad \text{  with } \quad \mathcal U_1(0,\cdot) \geq u_0(\cdot).
$$
 The exponential part of $\mathcal U_1$
is used only for $x\geq\xi_1(t)>c_1t$, where the reaction is $f_+$.
From Lemma \ref{Lem-super},   $f_+(t,s)\leq r_+(t)$, and
\eqref{eq-dissp}, we know that  $\mathcal{U}_1$ is a  generalized supersolution.  
Similarly to the previous upper estimates,  the comparison principle implies 
$$
 \overline c(u)\leq2\sqrt{R_+}.
$$

If $2\sqrt{R_+}<c_1<2\sqrt{R_-}$, set
$$
 \lambda_+:=\frac{c_1-\sqrt{c_1^2-4R_+}}2,
 \qquad \lambda_++\frac{R_+}{\lambda_+}=c_1.
$$
For $B$ sufficiently large, we define
$$
 \xi_2(t):=c_1t+\frac1{\lambda_+}
 \log\frac{B\theta_{r_+}(t)}M>c_1t  \qquad \forall t \ge 0,
$$
and
$$
 \mathcal U_2(t,x):=\min\{M,B\phi_{r_+,\lambda_+}(t,x)\} \qquad \text{ with } U_2(0,\cdot) \ge u_0(\cdot).
$$
One verify that $\mathcal{U}_2$ is a generalized supersolution. Hence, the comparison principle implies 
$$
 \lim_{t\to \infty} \sup_{x\geq ct}u(t,x)
 \leq  \lim_{t\to \infty}  B\max_{s\in[0,T]}\theta_{r_+}(s) 
 e^{-\lambda_+(c-c_1)t}=0 \qquad  \forall c>c_1.
$$
So we have $\overline c(u)\leq c_1$.

Finally, if $c_1\geq2\sqrt{R_-}$, the decreasing order gives
$$
 g(t,z,s)\leq f_-(t,s)\leq r_-(t).
$$
For $B$ large enough,
$$
 \mathcal U_3(t,x):=\min\{M,B\phi_{r_-,\sqrt{R_-}}(t,x)\}
$$
is a generalized supersolution. We have  $\overline c(u)\leq2\sqrt{R_-}$.  
The three estimates show
that
\begin{equation}\label{eq:dec-combined-upper}
 \overline c(u)\leq C_\downarrow(c_1;R_-,R_+).
\end{equation}

We next prove the lower estimates.  For all $z\in\R$ and $0\leq s\leq M$,
$$
 g(t,z,s)\geq f_+(t,s)\geq r_+(t)-Ks.
$$
Lemma \ref{Lem-lowerbound}, with $b=r_+$, therefore gives
\begin{equation}\label{eq:dec-global-lower-speed}
 \underline c(u)\geq2\sqrt{R_+}.
\end{equation}

It remains to use the lower bound
$g(t,x-c_1t,s)=f_-(t,s)\geq r_-(t)-Ks$ on $x<c_1t$.  Fix
$$
 0<q_0<q<\min\{c_1,2\sqrt{R_-}\}.
$$
Choose $l_0,l_q>0$ such that
$$
 l_0^{-2}<R_-,\qquad \frac{q^2}{4}+l_q^{-2}<R_-.
$$
For sufficiently small $\eta>0$, set
$$
 D_0=-1-l_0(\pi-\arcsin\eta),\qquad
 \Psi_0=\eta\psi_{r_-,0,0,l_0,D_0},\qquad
 P=\eta^2\theta_{r_-}.
$$
Define $h_q(y):=e^{-qy/2}\sin(y/l_q)$ for
$0\leq y\leq\pi l_q$, and let $y_q$ be the first solution of
$h_q(y)=\eta$ on the increasing part of $h_q$.  Put
$$
 \Psi_q=\eta\psi_{r_-,0,q,l_q,0},
 \qquad \xi_q(t)=qt+y_q.
$$
For $t_0$ sufficiently large,
$$
 -1<\xi_q(t)<qt+\pi l_q<c_1t,\qquad \forall  t\geq t_0.
$$
Define
$$
 \underline u(t,x):=
 \begin{cases}
 0,&x\leq D_0,\\
 \Psi_0(t,x),&D_0<x<-1,\\
 P(t),&-1\leq x\leq\xi_q(t),\\
 \Psi_q(t,x),&\xi_q(t)<x<qt+\pi l_q,\\
 0,&x\geq qt+\pi l_q
 \end{cases}
$$
The support of $\underline{u}$  is contained in $\{x<c_1t\}$ for $t\geq t_0$.  On this domain, the
reaction is $f_-$, and $f_-(t,s)\geq r_-(t)-Ks$. 
Let $\eta>0$ be sufficiently small. We can verify that the positive derivative gap  conditions are satisfied and $\underline{u}$ is a generalized subsolution.     Fix $t_0$. Decreasing $\eta$ such that $\underline u(t_0,\cdot)\leq u(t_0,\cdot)$, the comparison principle yields 
$$
 \liminf_{t\to\infty}\inf_{0\leq x\leq q_0t}u(t,x)
 \geq\eta^2\min_{[0,T]}\theta_{r_-}>0.
$$
By the arbitrariness of $q_0$, we have 
$$
 \underline c(u)\geq\min\{c_1,2\sqrt{R_-}\}.
$$
Together with \eqref{eq:dec-global-lower-speed}, this gives
\begin{align*}
 \underline c(u) \geq C_\downarrow(c_1;R_-,R_+).
\end{align*}
Combining this estimate with \eqref{eq:dec-combined-upper}, we complete the proof  of  
Theorem \ref{Thm-step}\textit{(ii)}.
\end{proof}

\section{Proof of Theorems \ref{Thm-inc} and \ref{Thm-dec}}
\label{Sec-general}
For the auxiliary monotone  nonlinear reactions of Section~\ref{Sec-step},
we first extend Theorem \ref{Thm-step} from linear shifts $c_1t$ to
shifts satisfying Assumption \ref{ASS-A}.

\subsection{Asymptotically linear shifts for stepwise profiles}
Let $g$ be defined in \eqref{stepg}, and let $A$ satisfy Assumption
\ref{ASS-A}.  Consider
\begin{equation}\label{eq-step2}
 \begin{cases}
 u_t=u_{xx}+u g(t,x-A(t),u),&t>0,\ x\in\R,\\
 u(0,x)=u_0(x),&x\in\R.
 \end{cases}
\end{equation}

\begin{proposition}\label{prop:step-shift}
Let $\csh$ be defined in \eqref{eq-A}, and let
 $u_0\in C(\R,[0,M])$ be nonzero and compactly supported.
The following statements hold:
\begin{itemize}
 \item[(i)] If $f_-(t,u)\leq f_+(t,u)$ for every
 $(t,u)\in[0,T]\times[0,M]$ and
 $0<R_-<R_+$, the spreading speed of \eqref{eq-step2} is
$$
 C_\uparrow(\csh;R_-,R_+).
$$
 \item[(ii)] If $f_-(t,u)\geq f_+(t,u)$ for every
 $(t,u)\in[0,T]\times[0,M]$ and $R_->R_+>0$, the spreading speed of
 \eqref{eq-step2} is
$$
 C_\downarrow(\csh;R_-,R_+).
$$
\end{itemize}
Here $C_\uparrow$ and $C_\downarrow$ are defined in
\eqref{eq:C-inc} and \eqref{eq:C-dec}, respectively.  The same
conclusion holds when $A(t)$ is replaced by $A(t)+D$ with fixed
$D\in\R$.
\end{proposition}

\begin{proof}
Fix $\delta\in(0,\csh)$.  By \eqref{eq-A}, there is
 $T_\delta>0$ such that
 \begin{equation*}
 (\csh-\delta/2)t\leq A(t)\leq(\csh+\delta/2)t,
  \qquad t\geq T_\delta.
 \end{equation*}
Define
\begin{equation*}
 \kappa_\delta:=\max\bigg\{0,\,\;
 \max_{0\leq t\leq T_\delta}
 \bigl((\csh-\delta)t-A(t)\bigr), \; \;\max_{0\leq t\leq T_\delta}
 \bigl(A(t)-(\csh+\delta)t\bigr)\bigg\}.
\end{equation*}
 By the
continuity of $A$,  the number $\kappa_\delta$ is nonnegative and finite.  
Thus
\begin{equation}\label{eq:A-sandwich}
  A_\delta^-(t):=(\csh-\delta)t-\kappa_\delta
  \leq A(t)\leq
  A_\delta^+(t):=(\csh+\delta)t+\kappa_\delta
\end{equation}
for every $t\geq0$.

If $f_-(t,u)\leq f_+(t,u)$, then
\begin{equation}\label{eq:inc-shift-reaction-order}
   g\bigl(t,x-A_\delta^+(t),u\bigr)
  \leq g \bigl(t,x-A(t),u\bigr)
  \leq g \bigl(t,x-A_\delta^-(t),u\bigr).
\end{equation}
Let $\underline u_\delta$ and
$\overline u_\delta$ be the solutions of
   \begin{equation*}
   \begin{cases}
   (\underline{u}_{\delta})_t= (\underline{u}_{\delta})_{xx}+ \underline{u}_{\delta} g\big( t, x - A_\delta^+(t),  \; \underline{u}_{\delta}\big), & t>0,\ x\in\R,\\
   \underline{u}_{\delta} (0,x)=u_0(x), & x\in\R,
   \end{cases}
   \end{equation*}
 and 
    \begin{equation*}
    \begin{cases}
    (\overline{u}_{\delta})_t= (\overline{u}_{\delta})_{xx}+ \overline{u}_{\delta} g\big( t, x - A_\delta^-(t), \;  \overline{u}_{\delta}\big), & t>0,\ x\in\R,\\
    \overline{u}_{\delta} (0,x)=u_0(x), & x\in\R,
    \end{cases}
    \end{equation*}
respectively.  By \eqref{eq:inc-shift-reaction-order} and comparison principle, one has 
$$
 \underline u_\delta\leq u\leq\overline u_\delta.
$$
The interfaces in the equations for $\underline u_\delta$ and
$\overline u_\delta$ are
$x=(\csh+\delta)t+\kappa_\delta$ and
$x=(\csh-\delta)t-\kappa_\delta$, respectively.  Theorem
\ref{Thm-step}   and comparison principle
 imply
$$
  C_\uparrow(\csh+\delta;R_-,R_+)
  \leq\underline c(u)\leq\overline c(u)
  \leq C_\uparrow(\csh-\delta;R_-,R_+),
$$
where $\underline c(u)$ and $\overline c(u)$ are the lower and upper
spreading  speeds  of \eqref{eq-step2}, respectively.
The continuity of $C_\uparrow$ gives
$$
 \underline c(u)=\overline c(u)=C_\uparrow(\csh;R_-,R_+).
$$

We next assume    that $f_- (t,u)\geq f_+(t,u)$.  The inequality
\eqref{eq:A-sandwich} gives
$$
  g \bigl(t,x-A_\delta^-(t),u\bigr)
  \leq g \bigl(t,x-A(t),u\bigr)
  \leq g \bigl(t,x-A_\delta^+(t),u\bigr).
$$
Let $\underline v_\delta$ and $\overline v_\delta$ denote the
solutions with initial value $u_0$ and shifts $A_\delta^-$ and
$A_\delta^+$, respectively.  Then
$$
 \underline v_\delta\leq u\leq\overline v_\delta.
$$
The two interfaces are $x=(\csh-\delta)t-\kappa_\delta$ and
$x=(\csh+\delta)t+\kappa_\delta$.  Theorem \ref{Thm-step}   yields
$$
  C_\downarrow(\csh-\delta;R_-,R_+)
  \leq\underline c(u)\leq\overline c(u)
  \leq C_\downarrow(\csh+\delta;R_-,R_+).
$$
Letting $\delta\to0$, we  have 
$$
 \underline c(u)=\overline c(u)=C_\downarrow(\csh;R_-,R_+).
$$
 
Finally,  for fixed $D\in \R$, one has
 $(A(t)+D)/t\to\csh$. 
 By the similar  constructions above to $A(t)+D$, we obtain the speed formula. 
The proof is complete. 
\end{proof}

\subsection{A comparison construction using only the linear growth}

We now return to the general equation \eqref{Pb}.   Let us put $f(t,z,u)$ between two stepwise reactions.

\begin{lemma}\label{lem:nonlinear-envelopes}
Suppose Assumptions~\ref{ASS-f} and \ref{ASS-limits} hold. Fix
\begin{equation}\label{eq:envelope-epsilon}
 0<\varepsilon<\min\{R_-,R_+,\sigma_0/2\}.
\end{equation}
Choose $L_\varepsilon>0$ such that
$|r(t,z)-r_\pm(t)|\leq\varepsilon$ when
$\pm z\geq L_\varepsilon$, uniformly in $t$. Under $(H_\uparrow)$,
set
\begin{align}
 a_\varepsilon^-(t,z)&=
 \begin{cases}r_-(t)-\varepsilon,&z<L_\varepsilon,\\
 r_+(t)-\varepsilon,&z\geq L_\varepsilon,\end{cases}\label{eq:linear-lower-inc}\\
 a_\varepsilon^+(t,z)&=
 \begin{cases}r_-(t)+\varepsilon,&z<-L_\varepsilon,\\
 r_+(t)+\varepsilon,&z\geq-L_\varepsilon.\end{cases}\label{eq:linear-upper-inc}
\end{align}
Under $(H_\downarrow)$, set instead
\begin{align}
 a_\varepsilon^-(t,z)&=
 \begin{cases}r_-(t)-\varepsilon,&z<-L_\varepsilon,\\
 r_+(t)-\varepsilon,&z\geq-L_\varepsilon,\end{cases}
 \label{eq:linear-lower-dec}\\
 a_\varepsilon^+(t,z)&=
 \begin{cases}r_-(t)+\varepsilon,&z<L_\varepsilon,\\
 r_+(t)+\varepsilon,&z\geq L_\varepsilon.\end{cases}
 \label{eq:linear-upper-dec}
\end{align}
In both cases, we have  $a_\varepsilon^-\leq r\leq a_\varepsilon^+$. Define
$$
 \begin{aligned}
 F_\varepsilon^-(t,z,u)&=a_\varepsilon^-(t,z)-Ku,\\
 F_\varepsilon^+(t,z,u)&=
 \min\{a_\varepsilon^+(t,z),\ L(M-u)-\sigma_0/2\}.
 \end{aligned}
$$
Then,   we have
$$
 F_\varepsilon^-(t,z,u)\leq f(t,z,u)\leq F_\varepsilon^+(t,z,u), \qquad  \forall (t,z,u) \in \R^2\times[0,M].
$$
Both auxiliary reactions  satisfy the
  hypotheses of Section~\ref{Sec-step}. Their linear growth
rates are $a_\varepsilon^-$ and $a_\varepsilon^+$,
respectively, and their nonlinear stepwise functions have the order required
in Theorem~\ref{Thm-step}.
\end{lemma}
\begin{proof}
Since  $r\geq r_-$ on
$z<L_\varepsilon$ and $r$ converges to $r_+$ on   $z\geq L_\varepsilon$,  \eqref{eq:linear-lower-inc} is well defined and $r\geq a_\varepsilon^-$.
For   \eqref{eq:linear-upper-inc} and $r\leq a_\varepsilon^+$, we  use the  convergence 
to $r_-$ on $z<-L_\varepsilon$ and $r\leq r_+$ elsewhere. 
Similarly, \eqref{eq:linear-lower-dec} and \eqref{eq:linear-upper-dec} are well defined and  $a_\varepsilon^- \le r \leq a_\varepsilon^+$.
 
The KPP assumption  gives $F_\varepsilon^- \leq f$.  The condition \eqref{eq-dissp} implies 
$$
 f(t,z,u)\leq f(t,z,M)+L(M-u)
 \leq L(M-u)-\sigma_0.
$$
Together with $f(t,z,u)\leq r(t,z)  \leq a_\varepsilon^+(t,z)$,  this implies
$f(t,z,u)\leq F_\varepsilon^+(t,z,u)$.  
From Assumption \ref{ASS-f} and \ref{ASS-limits}, we note that
$$r(t,x) \le \min\{   LM-\sigma_0, \; KM-\sigma_0\}  \text{  and }  r_{\pm} (t) \le \min\{   LM-\sigma_0, \; KM-\sigma_0\}.  $$
Combining with  \eqref{eq:envelope-epsilon}, one has 
$$
 a_\varepsilon^+(t,z)
 \leq LM-\sigma_0+\varepsilon<LM-\sigma_0/2.
$$
Note that  $F^{+}_\varepsilon (t,z,u)$ is Lipschitz continuous in $u$ with constant $L$ uniformly in $(t,z)$ and $F^{+}_\varepsilon (t,z,M) \leq -\sigma_0/2$.  
So  we have 
$ F_\varepsilon^+(t,z,0)=a_\varepsilon^+(t,z)$ and 
$$
 a_\varepsilon^+(t,z)-Lu
 \leq F_\varepsilon^+(t,z,u)\leq a_\varepsilon^+(t,z).
 $$
Similarly, we notice 
$$
 F_\varepsilon^-(t,z,M)
 =a_\varepsilon^-(t,z)-KM\leq-\sigma_0-\varepsilon.
$$
Both $F_\varepsilon^-$ and $F_\varepsilon^+$ satisfy the  hypotheses of Section~\ref{Sec-step}.
 Because the maps $a\mapsto a-Ku$ and $a\mapsto\min\{a,L(M-u)-\sigma_0/2\}$ are nondecreasing, under $(H_{\uparrow})$, we have 
 \begin{equation*}
 \begin{aligned}
 F^-_\varepsilon(t, -\infty, 0)    \leq F^-_\varepsilon(t, +\infty, 0),  \qquad  & F^+_\varepsilon(t, -\infty, 0)    \leq F^+_\varepsilon(t, +\infty, 0),    \\
\av{F^-_\varepsilon(\cdot, -\infty, 0) } = R_- - \varepsilon, \qquad  &  \av{F^-_\varepsilon(\cdot, +\infty, 0) } = R_+ - \varepsilon.
 \end{aligned}
 \end{equation*}
 Under $(H_{\downarrow})$, we obtain the reverse inequality.  The proof is complete.
 
\end{proof}

\begin{proof}[Proof of Theorems~\ref{Thm-inc} and \ref{Thm-dec}]
By  Lemma~\ref{lem:nonlinear-envelopes}, we consider
$$
 (u_\varepsilon^\pm)_t=(u_\varepsilon^\pm)_{xx}
 +u_\varepsilon^\pm F_\varepsilon^\pm
       (t,x-A(t),u_\varepsilon^\pm),\qquad
 u_\varepsilon^\pm(0,\cdot)=u_0.
$$
 The comparison principle implies 
$$
 0\leq u_\varepsilon^-\leq u\leq u_\varepsilon^+\leq M.
$$
Each auxiliary equation is a stepwise problem with jump at 
$A(t)+D_\varepsilon$, where $D_\varepsilon$ is either $L_\varepsilon$ or $-L_\varepsilon$. Applying Proposition~\ref{prop:step-shift} and Lemma~\ref{lem:nonlinear-envelopes},  we have   
\begin{equation*}
 \begin{split}
  C_\uparrow(\csh; R_{-} -\varepsilon,R_+ -\varepsilon)
  &\leq\underline c(u)\leq\overline c(u)\leq C_\uparrow (\csh;R_-+\varepsilon,R_++\varepsilon),\\
  C_\downarrow(\csh;R_{-} -\varepsilon,R_+ -\varepsilon)
  &\leq\underline c(u)\leq\overline c(u)\leq C_\downarrow (\csh;R_- +\varepsilon,R_+ +\varepsilon).
  \end{split}
\end{equation*}
Let  $\varepsilon\to 0$. By the continuity of $C_{\uparrow}$ and $C_{\downarrow}$, we complete the proof.

\end{proof}

\section{Extensions}\label{Sec-consequences}

\subsection{Spatially homogeneous time-periodic diffusion}
\begin{corollary}
Replace \eqref{Pb} by
$$
 u_t=d(t)u_{xx}+u f(t,x-A(t),u),
$$
where $d$ is continuous, $T$-periodic and
$0<d_{\min}\leq d(t)\leq d_{\max}<\infty$.
Suppose Assumptions  \ref{ASS-f}, \ref{ASS-limits} and \ref{ASS-A} hold,  and let
$u_0\in C(\R,[0,M])$ be nonzero and compactly supported. Denote $\overline d=\av d$.
The corresponding solution has rightward spreading speed
\begin{equation}\label{eq:periodic-diffusion-speed}
 \cstar=\overline d\,
 C_\diamond\left(\frac{\csh}{\overline d};
                 \frac{R_-}{\overline d},\frac{R_+}{\overline d}\right),
\end{equation}
with $\diamond=\uparrow$ or $\downarrow$ according to the order.  Particularly, 
in the increasing intermediate regime,  this reads
$$
 \cstar=\overline d\lambda+\frac{R_-}{\lambda},\qquad
 \lambda=\frac{\csh}{2\overline d}
         -\sqrt{\frac{R_+-R_-}{\overline d}},
$$
where
$$
 2\sqrt{\overline d R_+}<\csh<
 2\bigl(\sqrt{\overline d(R_+-R_-)}+\sqrt{\overline d R_-}\bigr).
$$
In the decreasing case,
$\cstar=\max\{2\sqrt{\overline d R_+},
                  \min\{\csh,2\sqrt{\overline d R_-}\}\}$.
\end{corollary}
\begin{proof}
Set $\tau(t)=\int_0^t d(s)\,ds$. Let $t=t(\tau)$ 
and  $v(\tau,x)=u(t(\tau),x)$. Then
$$
 v_\tau=v_{xx}+v\widetilde f(\tau,x-\widetilde A(\tau),v),\quad
 \widetilde f(\tau,z,v)=\frac{f(t(\tau),z,v)}{d(t(\tau))},\quad
 \widetilde A(\tau)=A(t(\tau)).
$$
Note that 
$$\widetilde r_\pm(\tau)= \frac{r_\pm(t(\tau))}{d(t(\tau))}  \quad \text{ with same period  } \mathcal T=\overline dT.
$$
This is because 
$t(\tau+\mathcal T)=t(\tau)+T$. Their means are
$$
 \frac1{\mathcal T}\int_0^{\mathcal T}\widetilde r_\pm(\tau)\,d\tau
 =\frac1{\overline dT}\int_0^T r_\pm(t)\,dt
 =\frac{R_\pm}{\overline d}.
$$
And  $\widetilde A(\tau)/\tau\to\csh/\overline d$ as $\tau \to\infty$. Apply the main
theorems  in the $\tau$ variable. We obtain
\eqref{eq:periodic-diffusion-speed}.  The proof is complete.
\end{proof}

\subsection{Shifts without a limiting slope}
\begin{theorem}
Let Assumptions  \ref{ASS-f} and \ref{ASS-limits} hold. Let $u_0\in C(\R, [0,M])$ be nonzero and compactly supported. Assume $A\in C([0,\infty))$ and
$$
 0<c_{\ell}:=\liminf_{t\to\infty}\frac{A(t)}t
 \leq c_h:=\limsup_{t\to\infty}\frac{A(t)}t<\infty.
$$
Then, one has 
\begin{align*}
 C_\uparrow(c_h;R_-,R_+)
 &\leq\underline c(u)\leq\overline c(u)
 \leq C_\uparrow(c_{\ell};R_-,R_+)
 &&\text{under }(H_\uparrow),\\
 C_\downarrow(c_{\ell};R_-,R_+)
 &\leq\underline c(u)\leq\overline c(u)
 \leq C_\downarrow(c_h;R_-,R_+)
 &&\text{under }(H_\downarrow).
\end{align*}
In particular, if the 
 function $C_\diamond(\cdot;R_-,R_+)$
is constant on $[c_{\ell},c_h]$, the spreading speed exists and equals
that constant, even when $A(t)/t$ does not converge.
\end{theorem}
\begin{proof}
For every $0<\delta<c_{\ell}$, continuity of $A$ and the definitions of $c_{\ell},c_h$  implies that there is  $\kappa_\delta >0$
such that
$$
 (c_{\ell}-\delta)t-\kappa_\delta
 \leq A(t)\leq(c_h+\delta)t+\kappa_\delta,
 \qquad  \forall  t\geq0. 
$$
From Lemma~\ref{lem:nonlinear-envelopes}, Theorem~\ref{Thm-step} and comparison principle, we have  
$$
 C_\uparrow(c_h+\delta;R_--\varepsilon,  R_+-\varepsilon) \leq\underline c(u)\leq\overline c(u) \leq C_\uparrow(c_{\ell}-\delta;R_-+\varepsilon,R_++\varepsilon).
$$
 Let $\varepsilon\to 0$ and 
$\delta\to0$.  The continuity of $C_\uparrow$ implies the conclusion under $(H_\uparrow)$.  Similarly, we can prove the conclusion under $(H_\downarrow)$. 
The last assertion follows when the lower and upper bounds coincide.  The proof is complete.
\end{proof}

\section{Discussion}
In this paper, we have established explicit rightward spreading speed formulas for Fisher-KPP equations in shifting environments where the limiting growth rates at the two far ends are time-periodic with positive averages. Our results reveal that  as long as the linear growth rate stays between its two spatial limits at each time, the long time invasion speed is completely determined by the two average growth rates and the asymptotic shift speed of the habitat. Temporary seasonal declines at the two ends, spatial nonmonotonicity and nonperiodic fluctuations inside the transition zone do not change the asymptotic spreading speed. Depending on which end is more favorable on average, an intermediate shift speed leads to either nonlocal pulling, where the population spreads faster than the locally linearly determined speed due to the better habitat ahead, or speed locking, where the invasion front moves at the same asymptotic speed as the environment.

 Our work relies on gluing explicit upper and lower solutions for an auxiliary stepwise problem and then comparing the general equation with shifted stepwise profiles. This is inspired by \cite[Chapter 6]{LamLou2022} for autonomous equations.
  Moreover, the main result extends naturally to asymptotically linear shifts and to spatially homogeneous, time-periodic diffusion through a direct change of the time variable.

Several natural questions remain open for future study. Our main assumptions require one limiting habitat to have a higher growth rate than the other at each time. If the growth rate inside the transition zone exceeds both spatial limits. For example, when a moving habitat offers a localized peak of resources, as studied in autonomous models \cite{GilettiGirardinMatano2024, LamNadinYu2025}, the transition zone itself may enhance  the speed of front and create new spreading regimes. It would be worthwhile to consider diffusion rates that vary with the moving frame, where a simple time rescaling no longer applies.

\section*{Acknowledgements}
 Z. Jin  is partially supported by the Strategic Priority Research  Program of CAS (No. XDB0900100), and NSFC (No. 12401225).  
Z. Wang  is  partially supported by the NSFC (Nos. 12301641 and 12471150).

\end{document}